# CONSTRUCTION OF $E(S^2)$-OPTIMAL SUPERSATURATED DESIGNS

By Dursun A. Bulutoglu and Ching-Shui Cheng[1]

*Air Force Institute of Technology, and Academia Sinica
and University of California, Berkeley*

Booth and Cox proposed the $E(s^2)$ criterion for constructing two-level supersaturated designs. Nguyen [*Technometrics* **38** (1996) 69–73] and Tang and Wu [*Canad. J. Statist* **25** (1997) 191–201] independently derived a lower bound for $E(s^2)$. This lower bound can be achieved only when $m$ is a multiple of $N-1$, where $m$ is the number of factors and $N$ is the run size. We present a method that uses difference families to construct designs that satisfy this lower bound. We also derive better lower bounds for the case where the Nguyen–Tang–Wu bound is not achievable. Our bounds cover more cases than a bound recently obtained by Butler, Mead, Eskridge and Gilmour [*J. R. Stat. Soc. Ser. B Stat. Methodol.* **63** (2001) 621–632]. New $E(s^2)$-optimal designs are obtained by using a computer to search for designs that achieve the improved bounds.

**1. Introduction.** In an experiment involving $m$ two-level factors, at least $m+1$ runs are required to estimate all the main effects. A design is called supersaturated if the run size is less than $m+1$. Under the assumption of effect sparsity that only a small number of factors are active, a supersaturated design can provide considerable cost saving in factor screening. Recently there have been quite a few articles on the analysis and construction of such designs. In particular, the $E(s^2)$ criterion proposed by Booth and Cox (1962) for constructing two-level supersaturated designs was studied by, for example, Lin (1993, 1995), Wu (1993), Nguyen (1996), Tang and Wu (1997), Cheng (1997), Li and Wu (1997), Butler, Mead, Eskridge and Gilmour (2001), Eskridge, Gilmour, Mead, Butler and Travnicek (2001) and Liu and Dean

Received October 2002; revised July 2003.

[1]Supported by NSF Grant DMS-00-71438 and National Security Agency Grant MDA904-00-1-0020.

(2002). This article contains further results on the construction of $E(s^2)$-optimal designs.

We represent an $N$-run supersaturated design for $m$ two-level factors by an $N \times m$ matrix $\mathbf{X}$ of 1's and $-1$'s. Each column corresponds to one factor and each row defines a factor-level combination. It is essential that no two columns of $\mathbf{X}$ are completely aliased, that is, there are no two columns $\mathbf{x}$ and $\mathbf{y}$ such that $\mathbf{x} = \mathbf{y}$ or $\mathbf{x} = -\mathbf{y}$. Throughout this article, we also assume that each column of $\mathbf{X}$ contains the same number of 1's and $-1$'s. Therefore, $N$ is even and the number of possible factors that can be accommodated is at most $\frac{1}{2}\binom{N}{N/2} = \binom{N-1}{N/2-1}$. Thus we have

$$(1.1) \qquad N - 1 < m \leq \binom{N-1}{N/2-1}.$$

The $E(s^2)$ criterion, in seeking a design as close to orthogonal as possible, minimizes $E(s^2) = \sum_{i<j} s_{ij}^2 / \binom{m}{2}$, where $s_{ij}$ is the $(i,j)$th entry of $\mathbf{X}^T\mathbf{X}$. Nguyen (1996) and Tang and Wu (1997) independently derived the following lower bound for any supersaturated design with $m$ factors and $N$ runs:

$$(1.2) \qquad E(s^2) \geq \frac{m - N + 1}{(m-1)(N-1)} N^2.$$

When $N \equiv 0 \pmod{4}$, this bound can be achieved only if $m$ is a multiple of $N - 1$; when $N \equiv 2 \pmod{4}$, $m$ needs to be an even multiple of $N - 1$. One question is whether this bound can be achieved for every multiple of $N - 1$ when $N \equiv 0 \pmod 4$ and every even multiple of $N - 1$ when $N \equiv 2 \pmod 4$. This appears to be very hard and the answer is yet unknown. One objective of this article is to provide some results in this direction.

In Section 2 we present a method for constructing designs that achieve the Nguyen–Tang–Wu bound. Section 3 contains the other main result of this article: improved lower bounds for $E(s^2)$ when the Nguyen–Tang–Wu bound is not achievable. It came to our attention that Butler, Mead, Eskridge and Gilmour (2001) also derived improved lower bounds for $E(s^2)$. Unlike our bounds, Butler, Mead, Eskridge and Gilmour's bounds do not apply to all cases; see the discussion in Section 3. We also report some new $E(s^2)$-optimal designs obtained by using a computer to search for designs that achieve the improved bounds. All the proofs are presented in Section 4.

Throughout this article, we use $\mathrm{GF}(s)$ to denote a finite field with $s$ elements. The multiplicative group that consists of the nonzero elements of $\mathrm{GF}(s)$ is cyclic, and a generator of this group is called a primitive element of the field. For each positive integer $q$, we denote the set $\{0, 1, \ldots, q-1\}$ of nonnegative integers less than $q$ by $\mathbb{Z}_q$. Multiplication and addition in $\mathbb{Z}_q$ are reduced modulo $q$ when necessary. For each subset $T = \{a_1, \ldots, a_t\}$ of $\mathbb{Z}_q$ and $b \in \mathbb{Z}_q$, the set $\{a_1 + b, \ldots, a_t + b\}$ is denoted by $T + b$.

**2. Construction of $E(s^2)$-optimal supersaturated designs via balanced incomplete block designs.** Let $m = q(N-1)$, where $q$ is a positive integer. Cheng (1997) showed that a supersaturated design that achieves the lower bound in (1.2) is equivalent to a balanced incomplete block design with $N-1$ treatments and $q(N-1)$ blocks of size $N/2-1$, abbreviated as $\text{BIBD}(N-1, q(N-1), N/2-1)$. Without loss of generality, we may assume that all the entries in the first row of a supersaturated design $\mathbf{X}$ are equal to 1. Let $\mathbf{Z}$ be obtained from $\mathbf{X}$ by deleting the first row. Then $\mathbf{Z} = [z_{ij}]_{(N-1) \times m}$ can be considered as the treatment-block incidence matrix of a binary incomplete block design with $N-1$ treatments and $m$ blocks of size $N/2-1$, where the $i$th treatment appears in the $j$th block if and only if $z_{ij} = 1$. Then $\mathbf{X}$ attains the lower bound in (1.2) if and only if $\mathbf{Z}$ is the treatment-block incidence matrix of a balanced incomplete block design. This equivalence of the existence of a $\text{BIBD}(N-1, q(N-1), N/2-1)$ and that of an $N \times q(N-1)$ $\mathbf{X}$ that attains the lower bound in (1.2) extends the well-known result that the existence of a $\text{BIBD}(N-1, N-1, N/2-1)$ is equivalent to that of an $N \times N$ Hadamard matrix. Note that for $\mathbf{X}$ to have no completely aliased columns, all the blocks of the corresponding BIBD must be distinct.

The largest $\text{BIBD}(N-1, q(N-1), N/2-1)$ with distinct blocks is the trivial one consisting of all the $\binom{N-1}{N/2-1}$ subsets of size $N/2-1$ of the $N-1$ treatments, which corresponds to the supersaturated design with the maximum number of factors given in (1.1). So the question raised in the Introduction about the existence of designs that achieve the Nguyen–Tang–Wu bound is equivalent to whether there exists a $\text{BIBD}(N-1, m, N/2-1)$ with distinct blocks for every $m$ satisfying (1.1) that is a multiple of $N-1$ when $N \equiv 0 \pmod{4}$ and an even multiple of $N-1$ when $N \equiv 2 \pmod{4}$. If there exists a Hadamard matrix of order $N$, then a $\text{BIBD}(N-1, q(N-1), N/2-1)$ exists for every $q$; a simple construction is to piece together $q$ $\text{BIBD}(N-1, N-1, N/2-1)$'s. The additional challenge in the supersaturated design construction is not to use the same block more than once. In fact, the construction of supersaturated designs by combining $q$ Hadamard matrices as proposed in Tang and Wu (1997) is equivalent to piecing together $q$ $\text{BIBD}(N-1, N-1, N/2-1)$'s. It is not clear how one can avoid duplicated blocks in their construction. Also, Tang and Wu's construction is applicable only when $N$ is a multiple of 4. Furthermore, what one needs is that the whole design is a BIBD; it does not have to be the union of $q$ $\text{BIBD}(N-1, N-1, N/2-1)$'s. Thus the construction based on BIBDs is more general and flexible. In this section we present a method of using difference families [Wilson (1972)] to construct $\text{BIBD}(N-1, q(N-1), N/2-1)$'s with distinct blocks. The readers are referred to Chapter VII of Beth, Jungnickel and Lenz (1999) for a discussion of the construction of BIBDs based on difference families.

The following is our first construction result.

THEOREM 2.1. *Suppose $N-1$ is an odd prime power, $q$ is an even divisor of $N-2$, $x$ is a primitive element of $\mathrm{GF}(N-1)$ and $T$ is a subset of $\mathbb{Z}_q$ of size $q/2$. Then the $q(N-1)$ sets $\{S_{r,a} : r = 0, \ldots, q-1, a \in \mathrm{GF}(N-1)\}$, where $S_{r,a} = \{x^{jq+i} + a : 0 \le j \le (N-2)/q - 1, i \in T+r\}$, form the blocks of a $\mathrm{BIBD}(N-1, q(N-1), N/2-1)$. Furthermore, if $(N-2)/q$ is odd and $U$ is a subset of $\mathbb{Z}_q$ of size $q/2$ such that $U^* = U + (q/2)$, where $U^*$ is the complement of $U$ in $\mathbb{Z}_q$, then the $q(N-1)/2$ sets $\{S_{r,a} : r \in U, a \in \mathrm{GF}(N-1)\}$ form the blocks of a $\mathrm{BIBD}(N-1, q(N-1)/2, N/2-1)$.*

Note that in the above construction, as defined earlier, the elements of $T+r$ are reduced modulo $q$ if necessary.

The BIBDs constructed in Theorem 2.1 may not have distinct blocks. For example, if $T + r = T$, then $S_{r,a} = S_{0,a}$. To construct BIBDs with distinct blocks, let $e$ be the smallest positive integer such that $T + e = T$. Then $e$, called the order of $T$, is a divisor of $q$. Let $\mathcal{T}_j = \{T + (j-1)e, \ldots, T + (j-1)e + e - 1\}$ for $1 \le j \le q/e$. Then $\mathcal{T}_1 = \cdots = \mathcal{T}_{q/e}$ and all the $e$ sets in each $\mathcal{T}_j$ are distinct. In this case, it can be seen that the design constructed in Theorem 2.1 consists of $q/e$ replications of an identical BIBD with distinct blocks as long as $q \ne N-2$. This leads to the following theorem.

THEOREM 2.2. *Suppose $N-1$ is an odd prime power, $q$ is an even divisor of $N-2$ with $q \ne N-2$, $x$ is a primitive element of $\mathrm{GF}(N-1)$ and $T$ is a subset of $\mathbb{Z}_q$ of size $q/2$. Let $e$ be the smallest positive integer such that $T + e = T$. Then the $e(N-1)$ sets $\{S_{r,a} : r = 0, \ldots, e-1, a \in \mathrm{GF}(N-1)\}$, where $S_{r,a} = \{x^{jq+i} + a : 0 \le j \le (N-2)/q - 1, i \in T+r\}$, are distinct and constitute a $\mathrm{BIBD}(N-1, e(N-1), N/2-1)$. Furthermore, if $(N-2)/q$ is odd and $U$ is a subset of size $e/2$ of $\{0, \ldots, e-1\}$ such that $U^* = U + (q/2)$, where $U^*$ is the complement of $U$ in $\{0, \ldots, e-1\}$ and the addition is reduced modulo $q$, then the $e(N-1)/2$ sets $\{S_{r,a} : r \in U, a \in \mathrm{GF}(N-1)\}$ constitute a BIBD with distinct blocks.*

The first design described in Theorem 2.2 is constructed by using $e$ initial (base) blocks $S_{0,0}, \ldots, S_{e-1,0}$, where $S_{r,0} = \{x^{jq+i} : 0 \le j \le (N-2)/q-1, i \in T+r\}$. The second design uses the $e/2$ initial blocks $S_{r,0}$, where $r \in U$. This construction is similar to Wilson's construction of balanced incomplete block designs as described in Theorem 5.2 of Beth, Jungnickel and Lenz [(1999), page 489]. However, Wilson did not consider the constraint of no repeated blocks. Also, a divisor $q$ of $N-2$ was used to construct BIBDs of block size $N/2 - 1$ in Theorems 2.1 and 2.2, while Wilson used it to construct BIBDs of block size $(N-2)/q$ or $(N-2)/q + 1$.

We note that if both $T$ and $U$ consist of all the integers $0 \le i \le q/2 - 1$, then $e = q$ and $U^* = U + (q/2)$. Thus a $\mathrm{BIBD}(N-1, q(N-1), N/2-1)$ with distinct blocks can be constructed for each even divisor $q$ of $N-2$ and if, in addition, $(N-2)/q$ is odd, then a BIBD$(N-1, q(N-1)/2, N/2-1)$ with distinct blocks also exists.

EXAMPLE 2.1. Let $N = 20$. Then $N - 2 = 18$ has two even divisors not equal to 18: 2 and 6. Since both $18/2$ and $18/6$ are odd, Theorem 2.2 can be used to construct BIBDs with 19, 38, 57 and 114 distinct blocks of size 9. The construction is based on the finite field GF(19), which is equivalent to $\mathbb{Z}_{19}$. Suppose we choose the primitive element 2. To construct a BIBD(19, 57, 9), we let the integer $q$ in Theorem 2.2 be 6, $T = U = \{0, 1, 2\}$ and use the three initial blocks $\{2^0, 2^1, 2^2, 2^6, 2^7, 2^8, 2^{12}, 2^{13}, 2^{14}\} = \{1, 2, 4, 7, 14, 9, 11, 3, 6\}$, $\{2^1, 2^2, 2^3, 2^7, 2^8, 2^9, 2^{13}, 2^{14}, 2^{15}\} = \{2, 4, 8, 14, 9, 18, 3, 6, 12\}$ and $\{2^2, 2^3, 2^4, 2^8, 2^9, 2^{10}, 2^{14}, 2^{15}, 2^{16}\} = \{4, 8, 16, 9, 18, 17, 6, 12, 5\}$. Adding the integers $0, 1, \ldots, 18$ (mod 19) to all the elements in the initial blocks produces a BIBD with 57 distinct blocks of size 9. Write down the $19 \times 57$ treatment-block incidence matrix in which the $(i, j)$th entry is equal to 1 if the $i$th treatment appears in the $j$th block and is equal to $-1$ otherwise. Then by adding a row of 1's to this treatment-block incidence matrix, one obtains an $E(s^2)$-optimal 20-run design for 57 factors.

Note that once the initial blocks are determined, except for the row of 1's, the other rows of the corresponding supersaturated design can be developed cyclically from an initial row. Eskridge, Gilmour, Mead, Butler and Travnicek (2001) and Liu and Dean (2002) also considered cyclic generation of $E(s^2)$-optimal and nearly optimal supersaturated designs.

The following result provides more flexibility in the construction of $E(s^2)$-optimal supersaturated designs.

THEOREM 2.3. *Let $q$ be an even divisor of $N - 2$ such that $q \neq N - 2$. Let $T$ and $T'$ be subsets of size $q/2$ of $\mathbb{Z}_q$ such that $T' \neq T + a$ for all elements $a$ of $\mathbb{Z}_q$. If $d_1$ and $d_2$ are BIBDs constructed by applying Theorem 2.2 to $T$ and $T'$, respectively, then $d_1$ and $d_2$ have no blocks in common; therefore their union is also a BIBD with distinct blocks.*

Let $\mathcal{F}$ be the set of all subsets of size $q/2$ of $\mathbb{Z}_q$. For any two such subsets $T$ and $T'$, we write $T \sim T'$ if there is an element $r \in \mathbb{Z}_q$ such that $T = T' + r$. Then clearly "$\sim$" is an equivalence relationship. Therefore, $\mathcal{F}$ is partitioned into disjoint equivalence classes. One can choose a set $T$ from each equivalence class to construct a BIBD according to the method of Theorem 2.2. The BIBDs constructed by using $T$'s from different equivalence classes have no blocks in common. Therefore, the union of these BIBDs is a BIBD with distinct blocks of size $N/2 - 1$, and can be used to construct $E(s^2)$-optimal supersaturated designs that attain bound (1.2). Note that the order of each set $T$ is equal to the size of the equivalence class that contains $T$.

6     D. A. BULUTOGLU AND C.-S. CHENGEXAMPLE 2.1 (continued). Again consider the case $N = 20$. More designs can be obtained by using Theorem 2.3. Take $q = 6$. The $\binom{6}{3} = 20$ subsets of size 3 of $\mathbb{Z}_6$ can be partitioned into four equivalence classes:

$$\{0,1,2\}, \{1,2,3\}, \{2,3,4\}, \{3,4,5\}, \{4,5,0\}, \{5,0,1\};$$
$$\{0,1,3\}, \{1,2,4\}, \{2,3,5\}, \{3,4,0\}, \{4,5,1\}, \{5,0,2\};$$
$$\{0,1,4\}, \{1,2,5\}, \{2,3,0\}, \{3,4,1\}, \{4,5,2\}, \{5,0,3\};$$
$$\{0,2,4\}, \{1,3,5\}.$$

Each set in one of the first three equivalence classes can be used to construct a BIBD(19, 57, 9) and a BIBD(19, 114, 9) with distinct blocks. By using one, two or all three of these equivalence classes, one can construct BIBD(19, 19$t$, 9)'s with distinct blocks for $t = 3, 6, 9, 12, 15$ and 18. The last equivalence class can be used to construct a BIBD(19, 19, 9) and a BIBD(19, 38, 9) with distinct blocks. Combining these with designs constructed from the first three equivalence classes, we obtain BIBD(19, 19$t$, 9)'s with distinct blocks and $E(s^2)$-optimal 20-run designs with 19$t$ factors, for $1 \leq t \leq 20$.

Note that in the above example, the designs constructed by using the last equivalence class are the same as those constructed by choosing $q = 2$.

EXAMPLE 2.2. Let $N = 18$. Then $N - 2 = 16$ has three even divisors not equal to 16: 2, 4 and 8. Let $q = 8$. It can be seen that the 70 subsets of size 4 of $\mathbb{Z}_8$ can be partitioned into 10 equivalence classes: 8 equivalence classes of size 8, 1 equivalence class of size 4 and 1 equivalence class of size 2. By applying Theorem 2.2, one can construct BIBD(17, 34$t$, 8)'s with no repeated blocks for $t = 1$ (using the equivalence class of size 2), 2 (using the equivalence class of size 4) and 4 (using each of the 8 equivalence classes of size 8). Combining the BIBDs constructed from different equivalence classes, one obtains BIBD(17, 34$t$, 8)'s with distinct blocks for all $t$'s such that $1 \leq t \leq 35$. Note that since 18 is not a multiple of 4, the number of blocks $b$ of a BIBD(17, $b$, 8) must be a multiple of 34.

EXAMPLE 2.3. For $N = 10$, applying Theorems 2.2 and 2.3 with $q = 4$, one can construct BIBD(9, 18$t$, 4)'s with no repeated blocks for $t = 1, 2, 3$. This is because the six subsets of size 2 of $\mathbb{Z}_4$ can be partitioned into two equivalence classes of sizes 4 and 2, respectively. In this case, the trivial BIBD has $\binom{9}{4} = 126$ blocks. Taking the complements of BIBD(9, 18$t$, 4)'s with no repeated blocks for $t = 1, 2$ and 3 in the trivial BIBD, we obtain BIBD(9, 18$t$, 4)'s with distinct blocks for $t = 4, 5$ and 6. This provides a complete solution of all BIBD(9, $b$, 4)'s with no repeated blocks, and thus all 10-run supersaturated designs that attain the Nguyen–Tang–Wu bound.

Sometimes one can also produce designs with distinct blocks by combining those constructed by using different even divisors of $N-2$. Let $T_i$ be a subset of size $q_i/2$ of $\mathbb{Z}_{q_i}$, $i=1,2$, where $q_1$ and $q_2$ are even divisors of $N-2$. Then by the same argument as in the proof of Theorem 2.2, one can show that the designs obtained by applying Theorem 2.2 to $T_1$ and $T_2$ have no blocks in common as long as all of the initial blocks of the two designs are different.

EXAMPLE 2.4. Let $N=14$. Then $N-2=12$ has three even divisors not equal to 12: 2, 4 and 6. As in Example 2.1, by choosing $q=6$, one can construct BIBD$(13, 26t, 6)$'s with distinct blocks for $t=1,3,4,6,7,9,10$. Note that since $12/6$ is even, fewer designs can be constructed here than in Example 2.1. For $q=4$, since $12/4$ is odd, a BIBD$(13,26,6)$ and a BIBD$(13,52,6)$ with distinct blocks can be constructed by applying Theorem 2.2 to the subset $T=\{0,1\}$ of $\mathbb{Z}_4$. It can be seen that these two designs have no blocks in common with any of those constructed by using $q=6$. It follows that one can construct BIBD$(13, 26t, 6)$'s with distinct blocks for $1 \le t \le 12$. As in Example 2.1, choosing $q=2$ does not produce new designs.

**3. Improved lower bounds for $E(s^2)$.** The following theorem presents some improved lower bounds for $E(s^2)$.

THEOREM 3.1. *Suppose $m$ is a positive integer such that $m > N-1$. Then there is a unique $q$ such that $-2N+2 < m - q(N-1) < 2N-2$ and $(m+q) \equiv 2 \pmod 4$. Let $g(q) = (m+q)^2 N - q^2 N^2 - mN^2$.*

1. *If $N \equiv 0 \pmod 4$, then*

$$E(s^2) \ge \begin{cases} \dfrac{g(q)+2N^2-4N}{m(m-1)}, & \text{when } |m-q(N-1)| < N-1, \\[2mm] \dfrac{g(q)-2N^2+4N+4N|m-q(N-1)|}{m(m-1)}, & \\ & \text{when } N-1 < |m-q(N-1)| \le \dfrac{3}{2}N-2, \\[2mm] \dfrac{g(q)+4N^2-4N}{m(m-1)}, & \text{when } |m-q(N-1)| > \dfrac{3}{2}N-2. \end{cases}$$

2. *If $N \equiv 2 \pmod 4$ and $q$ is even, then $E(s^2) \ge \max(h(q), 4)$, where*

$$h(q) = \begin{cases} \dfrac{g(q)+2N^2-4N+8}{m(m-1)}, & \text{when } |m-q(N-1)| < N-1, \\[2mm] \dfrac{g(q)-2N^2+20N+(4N-8)|m-q(N-1)|-24}{m(m-1)}, & \\ & \text{when } N-1 < |m-q(N-1)| \le \dfrac{3}{2}N-3, \\[2mm] \dfrac{g(q)+4N^2-4N}{m(m-1)}, & \text{when } |m-q(N-1)| > \dfrac{3}{2}N-3. \end{cases}$$

3. *If $N \equiv 2 \pmod{4}$ and $q$ is odd, then $E(s^2) \geq \max(h(q), 4)$, where*

$$h(q) = \begin{cases} \dfrac{g(q) + 2N^2 - 4N}{m(m-1)}, & \text{when } |m - q(N-1)| < N - 1, \\ \dfrac{g(q) - 2N^2 + 4N + 4N|m - q(N-1)|}{m(m-1)}, & \\ & \text{when } N - 1 < |m - q(N-1)| \leq \dfrac{3}{2}N - 1, \\ \dfrac{g(q) + 4N^2 - 12N + 8|m - q(N-1)| + 8}{m(m-1)}, & \\ & \text{when } |m - q(N-1)| > \dfrac{3}{2}N - 1. \end{cases}$$

[Butler, Mead, Eskridge and Gilmour](#) (2001) also derived some lower bounds for $E(s^2)$. Write $m$ as $m = q'(N-1) + r$, where $|r| < N/2$. Their result does not apply to the case where $N \equiv 2 \pmod{4}$ and $q'$ is odd, while our bounds apply to all cases. A numerical comparison suggests that their bounds agree with ours in the cases where they are applicable. Table 1 shows values of the Nguyen–Tang–Wu bound and our improved bound for the range $N \leq m \leq 2(N-1)$, where $N = 10, 12, 14$ and $16$.

We have used the computer software Gendex developed by [Nguyen](#) (1996) to search for $E(s^2)$-optimal designs. In many cases, we were able to find designs which achieve the improved bounds. Since citetr3 have reported $E(s^2)$-optimal 12- and 16-run designs, we list in Table 2 the new 10- and 14-run $E(s^2)$-optimal designs we have found.

**4. Proofs.**

Table 1
*Nguyen–Tang–Wu bound and the bound of Theorem* 3.1
*for $N \leq m \leq 2(N-1)$, $N = 10, 12, 14$ and $16$*

| $N$ | $m$ | Bound of Theorem 3.1 | Nguyen–Tang–Wu bound |
|---|---|---|---|
| 10 | 10 | 4 | 1.23456 |
|    | 11 | 4 | 2.22222 |
|    | 12 | 4 | 3.03030 |
|    | 13 | 4.61538 | 3.70370 |
|    | 14 | 5.05494 | 4.27350 |
|    | 15 | 5.52381 | 4.76190 |
|    | 16 | 5.86666 | 5.18518 |
|    | 17 | 5.88235 | 5.55555 |
|    | 18 | 5.88235 | 5.88235 |

TABLE 1
*Continued*

| $N$ | $m$ | Bound of Theorem 3.1 | Nguyen–Tang–Wu bound |
|---|---|---|---|
| 12 | 12 | 2.18181 | 1.19008 |
|    | 13 | 3.69230 | 2.18181 |
|    | 14 | 4.21978 | 3.02097 |
|    | 15 | 4.57142 | 3.74026 |
|    | 16 | 5.20000 | 4.36363 |
|    | 17 | 5.64705 | 4.90909 |
|    | 18 | 5.96078 | 5.39037 |
|    | 19 | 6.45614 | 5.81818 |
|    | 20 | 6.82105 | 6.20957 |
|    | 21 | 6.85714 | 6.54545 |
|    | 22 | 6.85714 | 6.85714 |
| 14 | 14 | 4       | 1.15976 |
|    | 15 | 4       | 2.15384 |
|    | 16 | 4       | 3.01538 |
|    | 17 | 4.94117 | 3.76923 |
|    | 18 | 5.67320 | 4.43438 |
|    | 19 | 6.05848 | 5.02564 |
|    | 20 | 6.35789 | 5.55465 |
|    | 21 | 6.66666 | 6.03076 |
|    | 22 | 6.90909 | 6.46153 |
|    | 23 | 7.41502 | 6.85314 |
|    | 24 | 7.82608 | 7.21070 |
|    | 25 | 7.84000 | 7.53846 |
|    | 26 | 7.84000 | 7.84000 |
|    | 27 | 8.38746 | 7.87692 |
|    | 28 | 8.80423 | 8.21728 |
|    | 29 | 8.82758 | 8.53333 |
|    | 30 | 8.82758 | 8.82758 |
| 16 | 16 | 2.13333 | 1.13777 |
|    | 17 | 3.76470 | 2.13333 |
|    | 18 | 4.18300 | 3.01176 |
|    | 19 | 4.49122 | 3.79259 |
|    | 20 | 5.38947 | 4.49122 |
|    | 21 | 6.09523 | 5.12000 |
|    | 22 | 6.64935 | 5.68888 |
|    | 23 | 7.08300 | 6.20606 |
|    | 24 | 7.42029 | 6.78261 |
|    | 25 | 7.68000 | 7.11111 |
|    | 26 | 7.87692 | 7.50933 |

PROOF OF THEOREM 2.1. Since $x$ is a primitive element of $\mathrm{GF}(N-1)$, we have $x^{N-2} = 1$, $x^{(N-2)/2} = -1$ and $1, x, x^2, \ldots, x^{N-3}$ are all distinct,

TABLE 2
*New* 10- *and* 14-*run* $E(s^2)$-*optimal designs*

$N = 10$, $m = 14$, $E(s^2) = 5.0549$

| −1 | 1 | 1 | 1 | −1 | 1 | −1 | −1 | −1 | −1 | 1 | −1 | 1 | −1 |
|---|---|---|---|---|---|---|---|---|---|---|---|---|---|
| 1 | 1 | −1 | 1 | 1 | −1 | −1 | −1 | 1 | 1 | 1 | −1 | −1 | 1 |
| −1 | −1 | −1 | 1 | 1 | −1 | −1 | 1 | −1 | −1 | −1 | 1 | 1 | 1 |
| −1 | −1 | 1 | −1 | −1 | −1 | 1 | −1 | −1 | 1 | −1 | −1 | −1 | 1 |
| 1 | 1 | −1 | −1 | −1 | −1 | 1 | 1 | −1 | 1 | −1 | −1 | 1 | −1 |
| 1 | −1 | 1 | −1 | 1 | 1 | −1 | 1 | −1 | −1 | 1 | −1 | −1 | −1 |
| −1 | 1 | −1 | −1 | −1 | 1 | 1 | 1 | 1 | −1 | 1 | 1 | −1 | 1 |
| 1 | 1 | 1 | 1 | 1 | −1 | 1 | −1 | 1 | −1 | −1 | 1 | −1 | −1 |
| 1 | −1 | 1 | 1 | −1 | 1 | 1 | 1 | 1 | 1 | 1 | 1 | 1 | 1 |
| −1 | −1 | −1 | −1 | 1 | 1 | −1 | −1 | 1 | 1 | −1 | 1 | 1 | −1 |

$N = 10$, $m = 15$, $E(s^2) = 5.5238$

| −1 | −1 | −1 | 1 | 1 | 1 | 1 | 1 | −1 | 1 | −1 | −1 | −1 | −1 | −1 |
|---|---|---|---|---|---|---|---|---|---|---|---|---|---|---|
| −1 | −1 | −1 | −1 | 1 | −1 | −1 | −1 | 1 | −1 | −1 | −1 | 1 | 1 | 1 |
| 1 | 1 | −1 | 1 | −1 | 1 | 1 | −1 | −1 | 1 | −1 | 1 | 1 | 1 | 1 |
| −1 | 1 | −1 | 1 | −1 | −1 | 1 | 1 | 1 | −1 | 1 | 1 | 1 | −1 | 1 |
| 1 | 1 | 1 | 1 | 1 | −1 | −1 | 1 | −1 | −1 | −1 | −1 | −1 | 1 | 1 |
| −1 | −1 | 1 | −1 | −1 | 1 | −1 | 1 | −1 | 1 | 1 | −1 | 1 | −1 | 1 |
| 1 | −1 | −1 | −1 | −1 | −1 | −1 | 1 | 1 | 1 | 1 | 1 | −1 | 1 | −1 |
| 1 | 1 | 1 | −1 | −1 | 1 | 1 | −1 | 1 | −1 | −1 | −1 | −1 | −1 | −1 |
| −1 | 1 | 1 | 1 | 1 | 1 | −1 | −1 | 1 | 1 | 1 | 1 | −1 | 1 | −1 |
| 1 | −1 | 1 | −1 | 1 | −1 | 1 | −1 | −1 | −1 | 1 | 1 | 1 | −1 | −1 |

$N = 14$, $m = 17$, $E(s^2) = 4.9412$

| −1 | 1 | 1 | −1 | −1 | −1 | 1 | 1 | 1 | −1 | −1 | 1 | −1 | −1 | −1 | −1 | 1 |
|---|---|---|---|---|---|---|---|---|---|---|---|---|---|---|---|---|
| 1 | −1 | −1 | 1 | 1 | −1 | 1 | −1 | 1 | −1 | 1 | 1 | −1 | −1 | 1 | 1 | −1 |
| −1 | −1 | 1 | 1 | 1 | 1 | 1 | 1 | 1 | −1 | −1 | −1 | 1 | 1 | 1 | 1 | 1 |
| 1 | −1 | −1 | −1 | 1 | −1 | 1 | 1 | −1 | 1 | 1 | −1 | 1 | 1 | −1 | −1 | 1 |
| 1 | −1 | −1 | −1 | −1 | 1 | −1 | 1 | 1 | 1 | −1 | −1 | −1 | −1 | −1 | 1 | −1 |
| 1 | 1 | −1 | −1 | 1 | −1 | −1 | −1 | 1 | 1 | −1 | 1 | 1 | 1 | 1 | 1 | 1 |
| −1 | −1 | 1 | 1 | −1 | −1 | −1 | −1 | 1 | 1 | 1 | −1 | 1 | −1 | 1 | 1 | 1 |
| −1 | 1 | −1 | −1 | −1 | 1 | 1 | 1 | −1 | 1 | 1 | 1 | 1 | 1 | 1 | 1 | −1 |
| −1 | 1 | −1 | 1 | 1 | 1 | −1 | 1 | 1 | 1 | 1 | −1 | −1 | −1 | 1 | −1 | 1 |
| −1 | 1 | 1 | −1 | 1 | −1 | −1 | −1 | −1 | 1 | −1 | 1 | −1 | −1 | 1 | 1 | −1 |
| −1 | −1 | 1 | −1 | 1 | 1 | 1 | −1 | −1 | 1 | −1 | 1 | −1 | −1 | 1 | −1 | −1 |
| 1 | 1 | 1 | 1 | −1 | −1 | −1 | 1 | −1 | −1 | −1 | −1 | −1 | 1 | 1 | −1 | −1 |
| 1 | 1 | 1 | 1 | −1 | 1 | 1 | −1 | 1 | 1 | 1 | −1 | 1 | 1 | −1 | −1 | −1 |
| 1 | −1 | −1 | 1 | −1 | 1 | −1 | −1 | −1 | −1 | 1 | 1 | −1 | −1 | −1 | −1 | 1 |

Table 2
*Continued*

$N = 14$, $m = 18$, $E(s^2) = 5.6732$

```
−1  1 −1  1  1 −1  1  1  1  1  1  1 −1 −1  1  1 −1
 1  1 −1 −1  1  1 −1 −1 −1 −1 −1  1 −1 −1 −1  1 −1 −1
 1  1 −1 −1 −1 −1  1  1  1 −1 −1  1 −1  1  1  1  1  1
 1  1  1  1  1  1 −1 −1  1  1 −1 −1 −1  1 −1 −1  1  1
 1 −1  1  1 −1  1 −1 −1  1  1  1  1  1  1  1  1  1 −1
−1 −1 −1  1 −1  1 −1  1  1 −1 −1 −1  1 −1 −1  1 −1  1
−1 −1  1 −1  1 −1 −1  1  1 −1  1 −1 −1 −1  1 −1  1 −1
−1 −1 −1  1  1 −1 −1 −1 −1 −1 −1  1  1  1  1 −1 −1  1
−1 −1  1 −1 −1  1  1  1 −1 −1 −1  1  1  1 −1 −1  1 −1
−1 −1 −1 −1  1  1  1 −1 −1  1  1 −1 −1  1  1  1  1  1
 1  1  1  1  1  1  1  1 −1 −1  1 −1  1 −1  1 −1 −1  1
 1 −1  1 −1 −1 −1  1 −1  1  1  1  1 −1 −1 −1 −1 −1  1
 1  1 −1 −1 −1 −1 −1  1 −1  1  1 −1  1  1 −1 −1 −1 −1
−1  1  1  1 −1 −1  1 −1 −1  1 −1 −1 −1 −1  1  1 −1 −1
```

$N = 14$, $m = 19$, $E(s^2) = 6.0585$

```
−1 −1  1 −1 −1 −1 −1 −1 −1  1 −1  1 −1  1  1  1  1  1  1
 1  1 −1 −1 −1 −1 −1  1  1  1 −1  1  1 −1  1 −1  1  1 −1
 1  1  1  1  1  1 −1 −1  1  1 −1 −1 −1 −1 −1 −1 −1  1  1
−1 −1 −1  1  1 −1  1 −1  1  1 −1  1  1 −1 −1  1 −1 −1 −1
−1 −1  1 −1 −1  1 −1  1 −1 −1  1  1 −1 −1 −1  1 −1 −1 −1
−1  1  1  1  1 −1  1  1 −1 −1 −1  1 −1  1 −1 −1  1  1 −1
 1  1  1 −1 −1 −1  1  1  1  1 −1  1  1 −1  1 −1 −1  1
−1 −1  1  1  1 −1 −1  1 −1  1  1 −1  1 −1  1 −1  1 −1  1
−1  1  1 −1  1  1  1 −1  1 −1  1 −1  1 −1  1  1  1  1 −1
 1 −1 −1  1 −1 −1 −1 −1  1 −1  1 −1 −1  1 −1 −1  1 −1 −1
 1  1 −1  1 −1  1  1  1 −1 −1 −1 −1 −1 −1  1  1  1 −1  1
 1  1 −1 −1  1  1  1 −1 −1  1  1  1 −1  1  1 −1 −1 −1 −1
 1 −1 −1  1  1  1 −1  1  1 −1  1  1  1  1 −1  1  1
−1 −1 −1 −1 −1  1  1 −1 −1 −1 −1 −1  1  1 −1 −1 −1  1  1
```

where 1 is the multiplicative identity. The multiset $\bigcup_{r=0}^{q-1}\bigcup_{u,v\in S_{r,0},\,u\neq v}\{u-v\}$ also can be expressed as

$$\bigcup_{\substack{i_1,i_2\in T \\ 0\leq j_1\leq (N-2)/q-1 \\ i_1\neq i_2 \text{ if } j_1=0}} \bigcup_{r=0}^{q-1} \bigcup_{j_2=0}^{(N-2)/q-1} \{x^{j_2q+i_2+r} - x^{(j_1+j_2)q+i_1+r}\},$$

where $i_1 + r$ and $i_2 + r$ are reduced modulo $q$ if necessary. For any fixed triple $(i_1, i_2, j_1)$, where $i_1, i_2 \in T$, $0 \leq j_1 \leq (N-2)/q - 1$ and $i_1 \neq i_2$ if $j_1 = 0$,

$$\bigcup_{r=0}^{q-1} \bigcup_{j_2=0}^{(N-2)/q-1} \{x^{j_2 q + i_2 + r} - x^{(j_1+j_2)q + i_1 + r}\}$$

$$= \bigcup_{r=0}^{q-1} \bigcup_{j_2=0}^{(N-2)/q-1} \{x^{j_2 q + i_2 + r}(1 - x^{j_1 q + [(i_1+r) - (i_2+r)]})\},$$

where again $i_1 + r$ and $i_2 + r$ are reduced modulo $q$ if necessary.

Now $\bigcup_{r=0}^{q-1} \bigcup_{j_2=0}^{(N-2)/q-1} \{x^{j_2 q + i_2 + r}\}$ covers every power of $x$ and hence every nonzero element of $\mathrm{GF}(N-1)$ exactly once. Since $1 - x^{j_1 q + [(i_1+r) - (i_2+r)]} \neq 0$ and there are $\frac{q}{2}(\frac{N-2}{2} - 1)$ triples $(i_1, i_2, j_1)$ such that $i_1, i_2 \in T$, $0 \leq j_1 \leq (N-2)/q - 1$ and $i_1 \neq i_2$ if $j_1 = 0$, it follows that $\bigcup_{r=0}^{q-1} \bigcup_{u,v \in S_{r,0},\, u \neq v} \{u - v\}$ covers each nonzero element of $\mathrm{GF}(N-1)$ $\frac{q}{2}(\frac{N-2}{2} - 1)$ times. Thus the $q$ sets in $\{S_{r,0} : 0 \leq r \leq q - 1\}$ are a difference family and the $q(N-1)$ sets in $\{S_{r,a} : 0 \leq r \leq q - 1,\, a \in \mathrm{GF}(N-1)\}$ constitute the blocks of a BIBD.

If $(N-2)/q$ is odd, then $(q + N - 2)/2 = \alpha q$ for some positive integer $\alpha$. Then since $x^{(N-2)/2} = -1$ and $U^* = U + (q/2)$,

$$\bigcup_{r \in U} \bigcup_{j_2=0}^{(N-2)/q-1} \{x^{j_2 q + i_2 + r} - x^{(j_1+j_2)q + i_1 + r}\}$$

$$= \bigcup_{r \in U^*} \bigcup_{j_2=0}^{(N-2)/q-1} \{x^{q/2} x^{(N-2)/2}(x^{(j_1+j_2)q + i_1 + r} - x^{j_2 q + i_2 + r})\}$$

$$= \bigcup_{r \in U^*} \bigcup_{j_2=0}^{(N-2)/q-1} x^{\alpha q} \cdot \{x^{(j_1+j_2)q + i_1 + r} - x^{j_2 q + i_2 + r}\}$$

$$= \bigcup_{r \in U^*} \bigcup_{j_2=0}^{(N-2)/q-1} \{x^{(j_1+j_2+\alpha)q + i_1 + r} - x^{(j_2+\alpha)q + i_2 + r}\}$$

$$= \bigcup_{r \in U^*} \bigcup_{j_2=0}^{(N-2)/q-1} \{x^{(j_1+j_2)q + i_1 + r} - x^{j_2 q + i_2 + r}\}.$$

This implies that $\bigcup_{r \in U} \bigcup_{u,v \in S_{r,0},\, u \neq v} \{u - v\} = \bigcup_{r \in U^*} \bigcup_{u,v \in S_{r,0},\, u \neq v} \{u - v\}$. Since $\bigcup_{r=0}^{q-1} \bigcup_{u,v \in S_{r,0},\, u \neq v} \{u - v\}$ covers each nonzero element of $\mathrm{GF}(N-1)$ $\frac{q}{2}(\frac{N-2}{2} - 1)$ times, each of $\bigcup_{r \in U} \bigcup_{u,v \in S_{r,0},\, u \neq v} \{u - v\}$ and $\bigcup_{r \in U^*} \bigcup_{u,v \in S_{r,0},\, u \neq v} \{u - v\}$ covers every nonzero element of $\mathrm{GF}(N-1)$ $\frac{q}{4}(\frac{N-2}{2} - 1)$ times. Thus the $q(N-1)/2$ sets in $\{S_{r,a} : r \in U, \, a \in \mathrm{GF}(N-1)\}$ constitute the blocks of a $\mathrm{BIBD}(N-1, q(N-1)/2, N/2-1)$. $\square$

PROOF OF THEOREM 2.2. Since $\mathcal{T}_1 = \cdots = \mathcal{T}_{q/e}$ and all the $e$ sets in each $\mathcal{T}_j$ are distinct, where $\mathcal{T}_j = \{T + (j-1)e, \ldots, T + (j-1)e + e - 1\}$, $1 \leq j \leq q/e$, the $e$ blocks $\{S_{r,0} : r = 0, \ldots, e-1\}$ are themselves a difference family and are all distinct. Assume that certain two blocks of the design $\{S_{r,a} : r = 0, \ldots, e-1, \, a \in \mathrm{GF}(N-1)\}$ are the same. Then there exist an integer $0 \leq r \leq e-1$ and $a \in \mathrm{GF}(N-1)$, either $r \neq 0$ or $a \neq 0$, such that the two sets $K_1 = \{x^{jq+i} + a : i \in T, 0 \leq j \leq (N-2)/q - 1\}$ and $K_2 = \{x^{jq+i} : i \in T + r, 0 \leq j \leq (N-2)/q - 1\}$ are the same. Then the sum of the elements of $K_1$ is equal to that of the elements of $K_2$. Since $q \neq N-2$, $1 - x^q \neq 0$; thus we have

$$\sum_{j=0}^{(N-2)/q-1} \sum_{i \in T} x^{jq+i} = \left[\sum_{j=0}^{(N-2)/q-1} x^{jq}\right]\left[\sum_{i \in T} x^i\right]$$
$$= \left[\frac{1-x^{N-2}}{1-x^q}\right]\left[\sum_{i \in T} x^i\right] = 0.$$

This implies that the sum of the elements of $K_1$ is equal to $\frac{1}{2}(N-2)a$ and the sum of the elements of $K_2$ is equal to 0. Therefore, $\frac{1}{2}(N-2)a = 0$ in $\mathrm{GF}(N-1)$, and hence $a = 0$. Now since $K_1 = K_2$, the two sets $\widehat{K}_1 = \{jq + i : i \in T, 0 \leq j \leq (N-2)/q - 1\}$ and $\widehat{K}_2 = \{jq + i : i \in T + r, 0 \leq j \leq (N-2)/q - 1\}$, with the elements being integers modulo $N-2$, must be equal. By the definition of $e$, this can happen only if $r = 0$, which is a contradiction. $\square$

Theorem 2.3 can be proved in the same way as Theorem 2.2.

PROOF OF THEOREM 3.1. First we state and prove a lemma.

LEMMA 1. *Let $\mathbf{y}$ be a $1 \times (N-1)$ vector with integer entries such that the first $p-1$ entries are congruent to $2 \pmod 4$, and the last $N-p$ entries are multiples of 4. Suppose $p \leq N/2$ and $m'$ is the sum of the entries of $\mathbf{y}$. If $N - 1 < |m'| < 2(N-1)$ and $\theta = (|m'| - 2p + 2)/4$ is an integer, then $0 < \theta < N - p$, and the sum of squares of the entries of $\mathbf{y}$ is minimized if and only if $\mathbf{y}$ has $p-1$ entries equal to $-2$, $\theta$ entries equal to $-4$ and $N - p - \theta$ entries equal to $0$ when $m' < 0$, or $p - 1$ entries equal to $2$, $\theta$ entries equal to $4$ and $N - p - \theta$ entries equal to $0$ when $m' > 0$.*

PROOF. It is enough to prove the case $m' < 0$. The other case follows by reversing the signs.

If $m' < 0$, then $\theta = (-m' - 2p + 2)/4 < (2N - 2 - 2p + 2)/4 < N - p$. Also, $|m'| > N - 1 \Rightarrow \theta > [N - 1 - 2p + 2]/4 \geq [N - 1 - N + 2]/4 > 0$. Therefore, $\theta$ is a positive integer less than $N - p$. If a vector has $p - 1$ entries equal to $-2$, $\theta$ entries equal to $-4$ and $N - p - \theta$ entries equal to 0, then the sum of all the entries is $-4\theta - 2(p - 1) = m'$. We shall show that a vector $\mathbf{y}^*$ with the smallest sum of squares of the entries among all vectors satisfying the conditions in the lemma must be of this form.

First we show that $\mathbf{y}^*$ cannot have positive entries. If not, let $y_i^*$ be a positive entry of $\mathbf{y}^*$. We claim that $\mathbf{y}^*$ has at least one entry, say $y_j^*$, such that $y_j^* \leq -4$. For otherwise, all the negative entries of $\mathbf{y}^*$ are greater than or equal to $-2$. Then since $p \leq N/2$, at most $N/2$ entries of $\mathbf{y}^*$ can be $-2$. It follows that $m' - y_i^* \geq (-2)N/2 = -N$, but $m' < -N + 1$ and $y_i^* \geq 2$ imply that $m' - y_i^* < -N - 1$, a contradiction. Therefore, there is at least one $y_j^*$ such that $y_j^* \leq -4$. Now replacing $y_i^*$ and $y_j^*$ with $y_i^* - 4$ and $y_j^* + 4$, respectively, keeps the sum of the entries of the vector unchanged, but $(y_i^* - 4)^2 + (y_j^* + 4)^2 < (y_i^*)^2 + (y_j^*)^2$ since $|y_i^* - 4| \leq |y_i^*|$ and $|y_j^* + 4| < |y_j^*|$. This means that $\mathbf{y}^*$ can be improved, contradicting the fact that $\mathbf{y}^*$ is optimal. Therefore, $\mathbf{y}^*$ cannot have positive entries.

Finally we show that $\mathbf{y}^*$ can have only entries from the set $\{-4, -2, 0\}$. Let $y_1^*, \ldots, y_{N-1}^*$ be the entries of $\mathbf{y}^*$. From the previous paragraph we know that all the entries of $\mathbf{y}^*$ are nonpositive. Thus $y_i^* \leq -2$ for all $1 \leq i \leq p - 1$ and $y_j^* \leq 0$ for all $p \leq j \leq N - 1$. We first show that $y_i^* = -2$ for all $1 \leq i \leq p - 1$. Compare $\mathbf{y}^*$ with the vector that has $p - 1$ entries equal to $-2$, $\theta$ entries equal to $-4$ and $N - p - \theta$ entries equal to 0. Since the sum of all the entries is a constant, we see that if there is a $y_i^* < -2$, $1 \leq i \leq p - 1$, then there must be at least one $y_j^* = 0$, where $p \leq j \leq N - 1$. Then $(y_i^* + 4)^2 + (y_j^* - 4)^2 = (y_i^*)^2 + 8y_i^* + 32 < (y_i^*)^2 + (y_j^*)^2$. The last inequality follows from $y_i^* \leq -6$. This again shows that $\mathbf{y}^*$ can be improved, which is not possible. Thus we must have $y_1^* = \cdots = y_{p-1}^* = -2$. Then the minimum of $\sum_{i=p}^{N-1}(y_i^*)^2$ subject to the constraint that all the $y_i^*$'s are multiples of 4 is attained when each $y_i^*$ is 0 or $-4$, $p \leq i \leq N - 1$. Since $\sum_{i=1}^{N-1} y_i^* = m'$, $\theta = (-m' - 2p + 2)/4$ entries must be equal to $-4$. $\square$

Now we are ready to prove the theorem. We denote the sum of squares of all the entries of a matrix $\mathbf{M}$ by $\mathrm{SS}(\mathbf{M})$. Then for a supersaturated design $\mathbf{X}$ with $m$ factors and $N$ runs, $E(s^2) = [\mathrm{SS}(\mathbf{X}^T\mathbf{X}) - mN^2]/[m(m-1)]$. A key fact used in Nguyen (1996) and Cheng (1997) is that $\mathrm{SS}(\mathbf{X}^T\mathbf{X}) = \mathrm{tr}[\mathbf{X}^T\mathbf{X}\mathbf{X}^T\mathbf{X}] = \mathrm{tr}[\mathbf{X}\mathbf{X}^T\mathbf{X}\mathbf{X}^T] = \mathrm{SS}(\mathbf{X}\mathbf{X}^T)$, and since each column of $\mathbf{X}$ has the same number of 1's and $-1$'s, $\mathbf{X}\mathbf{X}^T$ has zero row sums.

If $N \equiv 2 \pmod{4}$, then all the entries of $\mathbf{X}^T\mathbf{X}$ are congruent to 2 (mod 4). In particular, all the off-diagonal entries have absolute values at least 2. Therefore, we have the simple lower bound

(4.1) $\qquad N \equiv 2 \pmod{4} \implies E(s^2) \geq 4.$

Let $q_1$ be the integer such that $N - 1 < m - q_1(N - 1) < 2N - 2$. Since exactly one of $m + q_1$, $m + q_1 + 1$, $m + q_1 + 2$ and $m + q_1 + 3$ is congruent to 2 (mod 4), there is a unique $q$ such that $-2N + 2 < m - q(N - 1) < 2N - 2$ and $m + q \equiv 2$ (mod 4). Let $\mathbf{X}^*$ be obtained by adding $q$ columns of 1's to $\mathbf{X}$ and let $\mathbf{J}$ be the $N \times N$ matrix of 1's. Then

$$(4.2) \qquad \mathbf{X}^*(\mathbf{X}^*)^T = \mathbf{XX}^T + q\mathbf{J},$$

and so $\mathrm{SS}(\mathbf{X}^*(\mathbf{X}^*)^T) = \mathrm{SS}(\mathbf{XX}^T) + q^2\mathrm{SS}(\mathbf{J}) + 2q \cdot$ (the sum of all entries of $\mathbf{XX}^T) = \mathrm{SS}(\mathbf{XX}^T) + q^2N^2$. The last equality follows from the fact that $\mathbf{XX}^T$ has zero row sums. Thus

$$(4.3) \qquad E(s^2) = [\mathrm{SS}(\mathbf{X}^*(\mathbf{X}^*)^T) - q^2N^2 - mN^2]/[m(m-1)].$$

A lower bound for $E(s^2)$ can be obtained by bounding $\mathrm{SS}(\mathbf{X}^*(\mathbf{X}^*)^T)$.

Without loss of generality, assume that each of the first $p$ rows of $\mathbf{X}^*$ has an even number of entries equal to 1 and each of the last $N - p$ rows of $\mathbf{X}^*$ has an odd number of entries equal to 1. We can also assume that $p \leq N/2$, since if needed we can change the signs of all the entries in a certain column of $\mathbf{X}^*$. Then $\mathbf{X}^*(\mathbf{X}^*)^T$ has the form

$$(4.4) \qquad \begin{bmatrix} \mathbf{A} & \mathbf{C} \\ \mathbf{C}^T & \mathbf{B} \end{bmatrix},$$

where $\mathbf{A}$ is $p \times p$, all the entries of $\mathbf{C}$ are multiples of 4 and all the entries of $\mathbf{A}$ and $\mathbf{B}$ are congruent to 2 (mod 4). This follows from the fact that $m + q$, the number of columns of $\mathbf{X}^*(\mathbf{X}^*)^T$, is congruent to 2 (mod 4).

Since $\mathbf{XX}^T$ has zero row sums and its diagonal entries are equal to $m$, the sum of the off-diagonal entries in each of its rows is $-m$. Therefore, by (4.2) the sum of the off-diagonal entries in each row of $\mathbf{X}^*(\mathbf{X}^*)^T$ is $-m + q(N - 1)$, which is a multiple of 4 if and only if $q$ is odd and $N \equiv 2$ (mod 4), since $m + q \equiv 2$ (mod 4). Then since all the entries of $\mathbf{C}$ in (4.4) are multiples of 4, $2(p - 1)$ is a multiple of 4 if and only if $q$ is odd and $N \equiv 2$ (mod 4). It follows that

(4.5)    $p$ is odd if $q$ is odd and $N \equiv 2$ (mod 4); otherwise, $p$ is even.

Now we first consider the case where $|m - q(N - 1)| < N - 1$. Since all the entries of $\mathbf{A}$ and $\mathbf{B}$ are congruent to 2 (mod 4), where $\mathbf{A}$ and $\mathbf{B}$ are as in (4.4), they all have absolute values at least 2. So $\mathrm{SS}(\mathbf{X}^*(\mathbf{X}^*)^T)$ is at least $(m + q)^2N + F(p)$, where

$$F(p) \equiv 4p(p-1) + 4(N-p)(N-p-1) = 8p^2 - 8Np + 4N^2 - 4N.$$

Since $F'(p) = 16p - 8N$ has a zero at $p = N/2$ and $F(p)$ is a convex function of $p$, by (4.5), $F(p)$ is minimized at $p = N/2$ if $N$ is a multiple of 4 or if $N \equiv 2$ (mod 4) and $q$ is odd, and at $p = N/2 - 1$ when $N \equiv 2$ (mod 4) and $q$ is even. From these observations we can calculate a lower bound for $\mathrm{SS}(\mathbf{X}^*(\mathbf{X}^*)^T)$:

$$(4.6) \quad \mathrm{SS}(\mathbf{X}^*(\mathbf{X}^*)^T) \geq \begin{cases} (m+q)^2 N + 2N^2 - 4N + 8, \\ \qquad \text{if } N \equiv 2 \text{ (mod 4) and } q \text{ is even}, \\ (m+q)^2 N + 2N^2 - 4N, \\ \qquad \text{otherwise}. \end{cases}$$

The various lower bounds for the case $|m - q(N-1)| < N - 1$ as stated in the theorem can be obtained by combining (4.6) with (4.1) and (4.3).

Next we consider the case $|m - q(N-1)| > N - 1$. By the discussion in the paragraph preceding (4.5), the sum of the off-diagonal entries in each row of $\mathbf{X}^*(\mathbf{X}^*)^T$ is $-m + q(N-1)$ and $(|-m + q(N-1)| - 2p + 2)/4$ is an integer. By Lemma 4.1, with $m' = -m + q(N-1)$, the sum of squares of the off-diagonal entries of the first $p$ rows of $\mathbf{X}^*(\mathbf{X}^*)^T$ is minimized if in each of these rows, $p - 1$ entries have absolute values equal to 2, $(|-m + q(N-1)| - 2p + 2)/4$ entries have absolute values equal to 4, and the rest are equal to 0. Thus $\mathrm{SS}(\mathbf{X}^*(\mathbf{X}^*)^T)$ is at least $(m+q)^2 N + F(p)$, where

$$\begin{aligned} F(p) &\equiv 4p(p-1) + 4(N-p)(N-p-1) \\ &\quad + 2 \cdot (-4)^2 \cdot p(|-m + q(N-1)| - 2p + 2)/4 \\ &= -8p^2 - 8Np + 4N^2 - 4N + 8p(|-m + q(N-1)|) + 16p. \end{aligned}$$

Since $F'(p) = -8(2p + N - |-m + q(N-1)| - 2)$ has a zero at $p = [|-m + q(N-1)| - N + 2]/2$ and $F$ is a concave function of $p$, by (4.5), $F(p)$ is minimized at 0 or $N/2$ [when $N \equiv 0$ (mod 4)], 0 or $N/2 - 1$ [when $N \equiv 2$ (mod 4) and $q$ is even], and 1 or $N/2$ [when $N \equiv 2$ (mod 4) and $q$ is odd]. Here we have used the fact that $0 < (|-m + q(N-1)| - N + 2)/2 < N/2$.

When $N \equiv 0$ (mod 4), $F(p)$ is minimized at 0 if $(|-m + q(N-1)| - N + 2)/2 > N/4$, that is, if $|-m + q(N-1)| > \frac{3}{2}N - 2$; otherwise, it is minimized at $p = N/2$. Similarly, when $N \equiv 2$ (mod 4) and $q$ is even, $F(p)$ is minimized at 0 if $|-m + q(N-1)| > \frac{3}{2}N - 3$; otherwise, it is minimized at $p = N/2 - 1$. When $N \equiv 2$ (mod 4) and $q$ is odd, $F(p)$ is minimized at 1 if $|-m + q(N-1)| > \frac{3}{2}N - 1$; otherwise, it is minimized at $p = N/2$. Lower bounds for $\mathrm{SS}(\mathbf{X}^*(\mathbf{X}^*)^T)$ based on these observations together with (4.1) and (4.3) establish the various lower bounds for the case $|m - q(N-1)| > N - 1$ as stated in the theorem. $\square$

Air Force Institute of Technology  
AFIT/ENC  
Building 640 Room 151C  
2950 Hobson Way  
Wright–Patterson  
  Air Force Base, Ohio 45433-7765  
USA  
e-mail: dursun@math.uic.edu

Institute of Statistical Science  
Academia Sinica  
Taipei 115, Taiwan  
Republic of China  
e-mail: cheng@stat.sinica.edu.tw